\documentclass[12pt]{elsarticle}

\usepackage{amsmath}
\usepackage{amsfonts}
\usepackage{amssymb}
\usepackage{amsfonts}
\usepackage{amsthm}
\usepackage{amscd}
\usepackage{mathtools}
\usepackage{commath}
\usepackage{lmodern}
\usepackage{color}
\usepackage{a4wide}
\usepackage{bm} \usepackage{todonotes}

\usepackage[draft]{fixme}
\fxsetup{theme=color}
\usepackage{hyperref}

\usepackage[english]{babel}
\usepackage[IL2]{fontenc}
\usepackage[utf8]{inputenc}
\newtheorem{theorem}{Theorem}[section]
\newtheorem{lemma}[theorem]{Lemma}

\newtheorem{corollary}[theorem]{Corollary}
\newtheorem{problem}[theorem]{Problem}
\theoremstyle{definition}
\newtheorem{definition}[theorem]
{Definition}
\theoremstyle{remark}
\newtheorem{remark}[theorem]{\upshape\bfseries Remark}
\newtheorem{example}[theorem]{\upshape\bfseries Example}

\newcommand{\R}{\mathbb{R}}
\renewcommand{\H}{\mathbb{H}}
\newcommand*{\qi}{\mathbf{i}}
\newcommand*{\qj}{\mathbf{j}}
\newcommand*{\qk}{\mathbf{k}}
\newcommand*{\ci}{\mathrm{i}}
\newcommand*{\Cj}[1]{{#1}^\ast}
\newcommand{\circlearc}{0\smallfrown1}

\DeclareMathOperator{\residue}{res}

\begin{document}

\begin{frontmatter}
\journal{}

\title{Closed Bounded Rational Framing Motions}

\author[1]{Hans-Peter Schröcker}
\ead{hans-peter.schroecker@uibk.ac.at}
\address[1]{University of Innsbruck, Department of Basic Sciences in Engineering Sciences, Technikerstr.~13, 6020 Innsbruck, Austria}

\author[2]{Zbyněk Šír}
\ead{zbynek.sir@karlin.mff.cuni.cz}
\address[2]{Charles University, Faculty of Mathematics and Physic, Sokolovská 83, Prague 186 75, Czech Republic}

\begin{abstract}
  We present a method for constructing all bounded rational motions that frame a space curve $\mathbf{r}(t)$. This means that the motion guides an orthogonal frame along the curve such that one frame axis is in direction of the curve tangent. Existence of (bounded) framing motions is equivalent to $\mathbf{r}(t)$ being a (bounded) rational Pythagorean Hodograph curve. In contrast to previous constructions that rely on polynomial curves with smooth self-intersection, our motions and curves are infinitely differentiable. To this end, we develop the theory of Pythagorean hodograph curves parameterized over the projective line. We also provide a simple geometric necessary and sufficient condition on the spherical part of the motion, given by the homogeneous quaternionic preimage of the Pythagorean hodograph curve, that ensures the existence of a corresponding bounded, rational, and even regular framing motion. The translation part comes from the speed distribution, which must be a special positive rational function. This can in practice be ensured by semidefinite optimization methods. We illustrate our findings with a number of examples.
\end{abstract}


\end{frontmatter}

\section{Introduction}
\label{sec:introduction}

Motions of classical kinematics are typically generated by mechanisms with revolute joints and are thus naturally bounded and closed. The mechanical constraints imposed by links and joints can be modeled by algebraic equations and the same is true for point trajectories which are algebraic curves. A typical example is the motion and trajectories generated by a planar four-bar linkage. Only in exceptional cases do mechanisms generate trajectory curves that are compatible with standards of Computer-Aided Geometric Design that demand \emph{rationality}. However, a couple of years ago, it turned out that it is possible to go in the opposite direction, starting with a rational motion or rational trajectories and constructing mechanisms from them \cite{hegedus13,li16}. These possibilities bridge a gap between mechanism science and Computer-Aided Geometric Design and account for renewed interest in the design of rational motions.

Framing motions of space curves saw a similar history. They arise naturally in the differential geometry of space curves, but neither the motion of the classic Frenet-Serret frame, which has already been investigated by G.~Darboux, cf. \cite[p. 181]{farouki08}, nor the motions of the Bishop frame \cite{bishop75} are generically rational and the same is true for rotation-minimizing framing motions \cite{klok86} or minimal twist motions \cite{farouki18}. Quite a lot of research effort was dedicated to finding curves with natural rational framing motions or to construct rational framing motions on rational curves \cite{wagner97,maurer99,wang08,barton10,farouki16b,farouki19}. It was shown in \cite{kalkan22} that there is a strong relation between the theory of rational framing motions and the theory of \emph{Pythagorean Hodograph} (PH) curves. Indeed, for a rational framing motion, the trajectory of the origin must be a PH curve and the same quaternion polynomial is used for its construction and for the rotation part of the motion.

PH curves were first introduced in the seminal papers \cite{farouki90c} (planar case) and \cite{farouki94a} (spatial case). The vast majority of research papers focus on \emph{polynomial} PH curves that have the advantage of being constructed directly via the integration of the hodograph \cite{farouki08}. Among many various constructions for PH curves, two recent papers \cite{FaroukiClosedC1,FaroukiClosedC2} are particularly relevant to our research. They are devoted to the construction of closed and bounded polynomial PH curves with one exceptional point where only $C^1$ respectively $C^2$ continuity is achieved. Consequently, a piecewise rational framing motion along this closed trajectory is available.

The literature devoted to rational PH curves is much more sparse. Their construction requires quite a different approach since integration of a rational hodograph does not necessarily generate a rational curve. This problem was bypassed in \cite{Pottmann95} where all planar rational PH curves were constructed using a dual geometrical approach. A generalization to spatial rational PH curves appeared in \cite{FaroukiSir} and was later improved and exploited in \cite{FaroukiSir2,krajnc1,krajnc2,krajnc3}. A different approach of a more algebraic flavor was used in \cite{LEE2014689} to construct a particular kind of planar rational PH curves. All spatial rational PH curves (containing the planar ones) were algebraically constructed in \cite{kalkan22,schroecker23} via solving a system of linear equations and in \cite{SchroeckerSir:_optimal_interpolation} by imposing zero residue conditions on hodograph. The rational PH curves with rational arc-length were studied in \cite{FaroukiSakkalis2019,SCHROCKER2024128842}. 

Our research is motivated by the recent publications \cite{FaroukiClosedC1, FaroukiClosedC2}. We propose a method to improve the smoothness of the constructed closed motions from $C^1$ or $C^2$, respectively, to  $C^\infty$. This necessitates the use of \emph{rational} PH curves. The central result, Theorem~\ref{th:convex-hull}, is a simple geometric characterization of all spherical rational motions that can be extended to a bounded rational framing motion of a suitable closed and bounded rational PH curve.

In order to produce useful curves with rational framing motion, care should be taken that the curve is cusp-free. This requires a strictly positive speed function, which can be encoded by strict positivity of a certain polynomial factor and ensured by semidefinite optimization. We will demonstrate that this computation method goes well with the simultaneous imposition of zero residue conditions for the resulting hodograph.

The remainder of this paper is organized as follows. In Section~\ref{sec:preliminaries} we review the necessary theory of motion polynomials and PH curves. We also motivate the extension of the PH curves to the projective line which is then described in detail in Section~\ref{sec:closed}. There we also prove Theorem~\ref{th:convex-hull}, our main theoretical result. The computational aspects of the construction are then addressed in Section~\ref{sec:construction}. It is followed by  several examples in Section~\ref{sec:examples}.

\section{Pythagorean Hodograph Curves and Motion Polynomials}
\label{sec:preliminaries}

In this section, we will review the elementary theory of polynomial and rational Pythagorean Hodograph (PH) curves, which are traditionally considered over a compact interval. We will explore the necessary steps to achieve smooth, regular, and closed PH curves and we will demonstrate the strong relationship between PH curves and rational framing motions.

\subsection{PH Curves and Their Reparameterizations}

PH curves in any dimension are polynomial or rational parametric curves $\mathbf{r}(t)$, characterized by having a polynomial or rational speed function:
\begin{equation}\label{defPH}
  \Vert \mathbf{r}'(t) \Vert = \sigma(t).
\end{equation}
In our discussion, we will focus on spatial curves, which include planar curves as a special case.

A spatial PH curve can be constructed from a quaternion-valued polynomial $\mathcal{A}(t)$ and a real function $\lambda(t)$ using the formula
\begin{equation}\label{int0}
  \mathbf{r}(t) = \int \lambda(t) \mathcal{A}(t) \mathbf{i} \Cj{\mathcal{A}}(t) \dif t.
\end{equation}
Polynomial PH curves are derived from polynomial $\lambda(t)$ without any restrictions. To obtain rational curves, $\lambda(t)$ must be chosen to be rational and subject to zero-residue conditions (linear equations in the coefficients of the numerator of $\lambda$) to ensure the rationality of the integral \eqref{int0}, see e.g., \cite{schroecker23}.

\begin{figure}
  \centering
\includegraphics[]{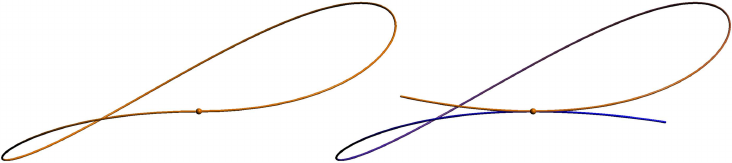}
  \caption{Closed curve constructed by methods of  \cite{FaroukiClosedC1} for $t\in [0,1]$ (left) and $t\in [-0.05,1.05]$ (right); varying color indicates actual course of the curve.}
\label{fig:1}
\end{figure}

The parameter $t$ is typically taken from a closed interval, most often $t \in [0,1]$. This is also the case in previous constructions of closed PH curves. In \cite{FaroukiClosedC1}, such closed curves are achieved by imposing the same position and tangent at the boundary parameter values $t=0$ and $t=1$. Consequently, the constructed closed periodic frames have one exceptional point and they exhibit $C^\infty$ smoothness at all points except one (the joint point), where only $C^1$ smoothness is achieved, c.f. Figure~\ref{fig:1}. A similar approach was used in order to obtain $C^2$ smoothness in~\cite{FaroukiClosedC2}.

Limiting the curve to a specific interval and using one privileged parameterization is quite natural in the case of polynomial PH curves. Indeed, in this case, the curve would necessarily be unbounded as the parameter value approaches infinity. Also, the only bijective re-parameterizations that preserve the polynomial nature of the curve are translations and scalings of the parameter domain. 

In the case of rational PH curves, however, it is natural to consider the group of rational linear parameter transformations
\begin{equation}\label{linrac}
  t = \psi(s) = \frac{as + b}{cs + d}.
\end{equation}

In this case, it is also pertinent to ask what happens as \( t \to \infty \). Indeed, the limit point of a rational curve can be finite and the curve can be perfectly smooth and regular in its neighbourhood. Moreover, the entire curve (considered over \( \mathbb{R} \)) can be bounded in some cases.

We will utilize all available parameter values of $t$, including the case where $t = \infty$; that is, we will parameterize the PH curves over the projective line. This is usual in algebraic geometry or rational kinematics but less so in Computer Aided Design. We will address issues involved in this construction properly in Section~\ref{sec:closed}. Here we just motivate it with three simple considerations.

Clearly, the PH property is independent of any particular interval for $t$. Indeed, due to the polynomial/rational nature of PH curves, if the equation \eqref{defPH} holds on a particular interval, then it holds over the entire line of real numbers~$\mathbb{R}$.

A second observation is that the rationality of the speed is not lost when the curve is rationally re-parameterized \cite{SirPlanarClassification}. Indeed, substituting $t = \psi(s)$, we get
\begin{equation}
  \label{eq:affine-reparam}
  \Bigl\Vert\od{}{s}\mathbf{r}(\psi(s))\Bigr\Vert =
  \Bigl\Vert \od{}{t}\mathbf{r}(t) \Bigr\Vert \od{}{s}\psi(s) =
  \sigma(\psi(s))\dot\psi(s),
\end{equation}
which is rational again, provided that $\psi(s)$ is rational. In particular for the reparameterization of the form \eqref{linrac}, we obtain
\begin{equation}\label{linracder}
  \dot{\psi}(s) = \frac{ad - bc}{(cs + d)^2}.
\end{equation}

Finally, if $\displaystyle \lim_{t \to \infty} \mathbf{r}(t)$ is a finite point, we can apply a suitable reparameterization for investigating the properties of the curve in its vicinity. As a simple but important example, consider the parameterization of the unit circle $\mathbf r_1(t)=\bigl(\frac{1-t^2}{t^2+1},\frac{2 t}{t^2+1}\bigr)$ with $\displaystyle \lim_{t \to \infty} \mathbf{r}_1(t)=(-1,0)$. Unlike the situation at Figure~\ref{fig:1} this point has the same regularity and smoothness as all other points on the circle. Taking simply $t=1/s$ we obtain the circle parameterization $\mathbf r_2(s) =\bigl( \frac{s^2-1}{s^2+1},\frac{2 s}{s^2+1}\bigr)$ with $\mathbf{r}_2(0)=(-1,0)$. Note that both circle parameterizations satisfy the PH property.

\subsection{Rational Framing Motions}

Formally speaking, a (rigid body) motion is a curve in $\operatorname{SE}(3)$. It is called rational if all trajectories are rational curves. Because of $\operatorname{SE}(3) = \operatorname{SO}(3) \ltimes \mathbb{R}^3$ it can be decomposed into a rotational part and a translation. The rotational part of a rational motion is given by a quaternionic polynomial $\mathcal{A}(t) \in \mathbb{H}[t]$ and maps the orthonormal basis $(\qi,\qj,\qk)$ to
\begin{equation}
  \label{eq:euler-rodriguez-frame}
  \mathbf{t}(t) \coloneqq \frac{\mathcal{A}(t)\qi\Cj{\mathcal{A}}(t)}{\mathcal{A}(t)\Cj{\mathcal{A}}(t)},\quad
  \mathbf{b}(t) \coloneqq \frac{\mathcal{A}(t)\qj\Cj{\mathcal{A}}(t)}{\mathcal{A}(t)\Cj{\mathcal{A}}(t)},\quad
  \mathbf{c}(t) \coloneqq \frac{\mathcal{A}(t)\qk\Cj{\mathcal{A}}(t)}{\mathcal{A}(t)\Cj{\mathcal{A}}(t)}.
\end{equation}
The translational part can be represented by a rational curve $\mathbf{r}(t)$. The point $\mathbf{p} = (x,y,z)$ is mapped to $\mathbf{r}(t) + x\mathbf{t}(t) + y\mathbf{b}(t) + z\mathbf{c}(t)$, the curve $\mathbf{r}(t)$ is the trajectory of the origin. The motion is said to be a \emph{framing motion} of $\mathbf{r}(t)$ if $\mathbf{t}(t)$ is parallel to the derivative vector $\mathbf{r}'(t)$, which implies that $\mathbf{r}(t)$ is a PH curve \cite{kalkan22}.

The construction of a PH curve $\mathbf{r}(t)$ via \eqref{int0} by integrating its hodograph $\lambda(t)\mathcal{A}(t)\qi\Cj{\mathcal{A}}(t)$ creates a framing motion whose rotational part is given by the polynomial $\mathcal{A}(t) \in \mathbb{H}[t]$ and the translation is $\mathbf{r}(t)$.  We also say that $\mathbf{r}(t)$ extends the spherical rational rotation motion $\mathcal{A}(t)$ to a rigid body motion. Note that this extension requires a choice of $\lambda(t)$.

The frame given by \eqref{eq:euler-rodriguez-frame} is usually called the \emph{Euler-Rodriguez frame} \cite{choi02}. The associated rigid body motion is called \emph{Euler-Rodriguez motion}. It is actually not unique as it can be composed from the right with a rational rotation around the tangent vector $\mathbf{t}(t)$. Algebraically, this can be realized by 
\begin{enumerate}
\item replacing the quaternionic $\mathcal{A}(t)$ by $\mathcal{A}(t)\mathcal{R}(t)$ where $\mathcal{R}(t) \in \mathbb{A}[t]$ is a quaternionic polynomial with coefficients in the algebra $\mathbb{A} \cong \mathbb{C}$ generated by $1$ and $\qi$ and
\item replacing the rational function $\lambda(t)$ by $\lambda(t)/(\mathcal{R}(t)\Cj{\mathcal{R}}(t))$.
\end{enumerate}
Indeed, because $\mathcal{R}(t)$ and $\qi$ commute, we have
\[
    \frac{\lambda(t)}{\mathcal{R}(t)\Cj{\mathcal{R}}(t)}
    (\mathcal{A}(t)\mathcal{R}(t))
    \qi
    \Cj{(\mathcal{A}(t)\mathcal{R}(t))}
    =
    \frac{\lambda(t)}{\mathcal{R}(t)\Cj{\mathcal{R}}(t)}
    \mathcal{R}(t)\Cj{\mathcal{R}(t)} \mathcal{A}(t)
    \qi
    \Cj{\mathcal{A}}(t)
    =
    \lambda(t)
    \mathcal{A}(t)
    \qi
    \Cj{\mathcal{A}}(t).
\]
In the light of this, it is probably more appropriate to speak of Euler-Rodriguez motions in case of rational PH curves. There is, however, a distinguished Euler-Rodriguez motion of the lowest possible degree from which all other Euler-Rodriguez motions can be obtained. Call the polynomial $\mathcal{A}(t) \in \H[t]$ \emph{$\qi$-reduced} if it has no right factor $\mathcal{R}(t) \in \mathbb{A}[t]$. Any PH curve $\mathbf{r}(t)$ can be generated by an $\qi$-reduced quaternionic polynomial $\mathcal{A}(t)$ and all rational framing motions of $\mathbf{r}$ are obtained by applying the replacements 1. and 2. above, c.f.~\cite{kalkan22}. The $\qi$-reduced polynomial $\mathcal{A}(t)$ is unique up to right multiplication with unit quaternions in $\mathbb{A}$ which change the frames along the curve by a fixed rotation about the tangent vector but not the framing motion itself.

To sum up, our construction of framing motions is based on the extension of rational rotational motion $\mathcal{A}(t)$ via a PH curve $\mathbf{r}(t)$ (the trajectory of the origin) to a rigid body motion. By \eqref{int0} this amounts to computing a suitable rational function $\lambda(t)$. It also necessitates to extend the parameter range from some compact interval to the projective line $\mathbb{P}_1$ or, equivalently, to the circle $\mathbb{S}_1$. We will discuss this in detail in Section~\ref{sec:closed}. There we also address the important issues of boundedness and regularity of~$\mathbf{r}(t)$.

\section{PH Curves Over the Projective Line}
\label{sec:closed}

In practice, bounded rational PH curves will often be constructed using parameterizations and functions over $\mathbb{R}$, see Section~\ref{sec:construction}. This section provides the theoretical foundations for our computations there. We feel that it is more suitable to
define all the objects over the real projective line $\mathbb{P}_1$ (or equivalently over the unit circle $\mathbb{S}_1$). This way the equal status of all points is obvious and the choice of the point corresponding to the infinite value of the parameter range is a purely formal aspect of the representation.

\begin{definition}
  A regular rational map $\mathbf R\colon \mathbb P_1\to \mathbb R^3$ is called \emph{regular bounded rational curve}. Throughout this paper it will also be understood as a map $\mathbf R\colon \mathbb S_1\to \mathbb R^3$ from the unit circle birationally identified with the projective line via
  \begin{equation}\label{circle}
    (u:v) \in \mathbb{P}_1 \longleftrightarrow \Bigl(\frac{u^2-v^2}{u^2+v^2},\frac{2uv}{u^2+v^2}\Bigr) \in \mathbb{S}_1
  \end{equation}
  and oriented by its unit tangent field \[\tau (u,v)=\Bigl(\frac{-2uv}{u^2+v^2},\frac{u^2-v^2}{u^2+v^2}\Bigr).\]
\end{definition}

\begin{remark}
  Comparing to standard rational curves in the sense of algebraic geometry we require all the points of the curve to be finite (no points at infinity). Moreover, we also wish to avoid cusps, hence the regularity assumption. Obviously $\mathbf R$ can be composed with any rational automorphism of $\mathbb P_1$ providing a different parameterization. Any embedding $\mathbb R \hookrightarrow \mathbb P_1$ provides a smooth parametric curve in the sense of differential geometry.

  Introduction of the reference embedding \eqref{circle} is not quite necessary but it provides the orientation and, more importantly, it allows a definition of the curve speed as a function. Without a particular embedding, the speed function would correspond rather to a finite measure or Riemannian metric defined over~$\mathbb P_1$.
\end{remark}

\begin{definition}\label{projPH}
  For a given regular bounded rational curve $\mathbf R\colon \mathbb S_1\to \mathbb R^3$ we denote by $\dif \mathbf R(u,v)$ its tangent map at the point $(u,v)$. We define its speed function
  \begin{equation*}
    L(u,v) \coloneqq \Vert \dif \mathbf R(u,v)(\tau(u,v)) \Vert
  \end{equation*}
  and its unit tangent vector
  \begin{equation*}
    \mathbf T(u,v) \coloneqq \frac{1}{L(u,v)} \dif \mathbf R(u,v)(\tau(u,v)).
  \end{equation*}
  We say that  $\mathbf R(u,v)$ is a \emph{regular bounded rational PH curve} if $L(u,v)$ is a rational function and consequently $\mathbf T(u,v)$ is a rational unit vector field.
\end{definition}

\begin{remark}
  Being Pythagorean hodograph in the sense of Definition \ref{projPH} does not depend on the reparameterization via the rational automorphisms of~$\mathbb P_1$. This can easily be seen similar to Equation~\eqref{eq:affine-reparam}.
\end{remark}

\begin{example}
  \label{ex:1}
  Let us consider the mapping
  \begin{equation*}
    \mathbf R(u,v)=\frac{1}{15 (u^2 + v^2)^5}
    \begin{pmatrix}
      14 u^{10}+130 u^8 v^2+110 u^6 v^4+30 u^4 v^6-60 u^2 v^8\\
      60 u^9 v+60 u^7 v^3+96 u^5 v^5+60 u^3 v^7+60 u v^9\\
      -120 u^8 v^2-300 u^6 v^4-300 u^4 v^6-120 u^2 v^8
    \end{pmatrix}.
  \end{equation*}
  It is a bounded regular rational curve (Figure~\ref{fig:2}).
\begin{figure}
    \centering
\includegraphics{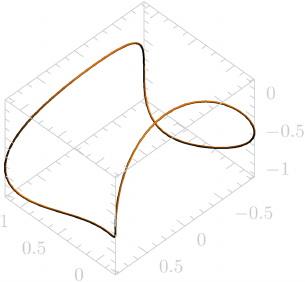}
    \caption{Regular, closed rational PH curve of Example~\ref{ex:1}.}
    \label{fig:2}
  \end{figure}
By a direct computation we get
  \begin{equation*}
    \dif \mathbf R(u,v)(\tau(u,v)) =
    \frac{1}{(u^2 + v^2)^5}
    \begin{pmatrix}
      4 u^9 v-20 u^7 v^3-16 u^5 v^5-20 u^3 v^7+4 u v^9\\
      2u^{10}-12 u^8 v^2+2 u^6 v^4-2 u^4 v^6+12 u^2 v^8-2 v^{10}\\
      -8u^9 v-8 u^7 v^3+8 u^3 v^7+8 u v^9
    \end{pmatrix}
  \end{equation*}
  yielding
  \begin{equation*}
    L(u,v)=\frac{2 (u^4+u^2 v^2+v^4)}{(u^2+v^2)^2},
  \end{equation*}
  whence $\mathbf{R}(u,v)$ is PH, and
  \begin{equation*}
    \mathbf T(u,v)=\frac{1}{(u^2 + v^2)^3}
    \begin{pmatrix}
      2 u^5 v-12 u^3 v^3+2 u v^5\\
      u^6-7 u^4 v^2+7 u^2 v^4-v^6\\
      4u v^5-4 u^5 v
    \end{pmatrix}.
  \end{equation*}
  Note that the speed function $L(u,v)$ is positive and bounded due to the regularity and boundness of the curve $\mathbf R$. This in particular means that it is homogeneous of degree $0$ and has no real roots in its denominator.
\end{example}

Similarly to the case of usual curves, the unit tangent vector field and the speed distribution fully define the curve (up to a translation). This observation is formalized in the following corollary.

\begin{corollary}
  \label{cor:zeroInt}
  The closed rational curve can be recovered from the functions $L(u,v)$ and $\mathbf T(u,v)$ via a line integral. More precisely, denoting by $\circlearc$ the oriented circle arc between the fixed point $(u_0,v_0)\in \mathbb S_1$ and some point $(u_1,v_1)\in \mathbb S_1$ we have
  \begin{equation}\label{int}
    \mathbf R(u_1, v_1) = \mathbf R(u_0, v_0) + \int \limits_{0\smallfrown 1}L(u,v)\mathbf T(u,v) \dif\mathbb S_1
  \end{equation}
  providing in particular
  \begin{equation}\label{zeroInt}
    \int_{\mathbb S_1}L(u,v)\mathbf T(u,v) \dif\mathbb S_1=\mathbf 0.
  \end{equation}
\end{corollary}

Because $\mathbf T(u,v)$ is a rational unit vector field, it can be seen as a trajectory under a spherical motion, which can always be expressed using a homogeneous quaternion polynomial $\mathcal A(u,v)$. We thus obtain a method for constructing all regular bounded rational PH curves from two inputs: homogeneous quaternion polynomial and degree $0$ homogeneous speed function. We formulate this fact in Theorem \ref{th:all-CRRPH-curves1} and complete it with the formulation of the practical problem of finding compatible $\mathcal A$ and $L$ in Problem \ref{problem}.

\begin{theorem}
  \label{th:all-CRRPH-curves1}
  All bounded rational PH curves can be obtained via formula \eqref{int} from a homogeneous quaternion polynomial $\mathcal A(u,v)$ and a positive bounded rational function $L(u,v)$ by setting
  \begin{equation}
    \label{eq:T}
    \mathbf T(u,v)=\frac{\mathcal A(u,v) \qi \mathcal A^*(u,v)}{\mathcal A(u,v) \mathcal A^*(u,v)}.
  \end{equation}
\end{theorem}

\begin{problem}\label{problem}
  For a chosen homogeneous quaternion polynomial $\mathcal A(u,v)$ find all the bounded strictly positive rational functions $L(u,v)$ so that \eqref{int} is rational.
\end{problem}

One of the key observations for solving Problem~\ref{problem} is that -- under some assumptions -- \eqref{zeroInt} implies the rationality of \eqref{int}. This will be shown with the help of the following

\begin{lemma}\label{zeroRac}
  Let $F(u,v)$ be a bounded rational function on the unit circle with only one pair of conjugate complex roots in its denominator and let its numerator be divisible by $(u^2+v^2)$. Then its integral
  \begin{equation}\label{int4}
    G(u_1,v_1) \coloneqq \int \limits_{0\smallfrown 1}F(u,v)\dif\mathbb S_1
  \end{equation}
  is rational if and only if
  \begin{equation}\label{int5}
    \int_{\mathbb S_1}F(u,v)\dif\mathbb S_1=\mathbf 0.
  \end{equation}
\end{lemma}

\begin{proof}
 Let us substitute $u=1$ and $v=t$, thus obtaining the following parameterization of the unit circle \eqref{circle}
  \begin{equation}
    \label{dS1}
    \mathbf c(t)=\Bigl( \frac{1-t^2}{1+t^2},\frac{2t}{1+t^2} \Bigr)
    \quad \text{and providing} \quad
    \dif\mathbb S_1=\Vert\mathbf c'(t)\Vert=\frac{2}{1+t^2}\dif t.
  \end{equation}
  By the lemma's assumption, the de-homogenization $F(1,t)$ of the function $F(u,v)$ is divisible by $1+t^2$ and it is of degree $0$ (due to the boundedness). Setting $f(t)=2F(1,t)/(1+t^2)$ we obtain a rational function of degree $-2$ with one pair of complex conjugate roots of its denominator. The integrals \eqref{int4}, \eqref{int5} become
  \begin{equation}
    \label{intvanish}
    G(1,t) = \int f(t)\dif t
    \quad \text{and} \quad
    \int_{\mathbb S_1}F(u,v)\dif\mathbb S_1=\int \limits_{-\infty}^\infty f(t)\dif t=\mathbf 0,
  \end{equation}
  respectively. The equivalence stated in the lemma now follows from the Weierstrass residue theorem due to the fact that the integrand has only one pole in the upper complex half-plane and its degree is~$-2$.
\end{proof}

A further technical lemma will be needed:

\begin{lemma}\label{lemBanach}
  Let $X$ be a Banach space, $Y \subset X$ a dense subspace and $b\colon X \to \mathbb R^n$ a linear map. Then $(\ker b)\cap Y$ is dense in $(\ker b)$.
\end{lemma}
\begin{proof}
  Suppose first that $n=1$ and let $B(x,\epsilon)$ denotes the open ball of radius $\epsilon$ centered at $x\in (\ker b)$. By density of $Y$ there exist points $y_1$, $ y_2\in Y \cap B(x,\epsilon)$ satisfying $b(y_1)>0$ and $b(y_2)<0$. The intersection of the segment $y_1y_2$ with $(\ker b)$ lies in both, $B(x,\epsilon)$ and $(\ker b)$, showing that $(\ker b)\cap Y$ is dense in $(\ker b)$. The statement for $n>1$ can now be proven by applying the case $n=1$ to each of the components of~$b$.
\end{proof}

We are now in a position to provide a simple geometric criterion for Problem~\ref{problem} to have solutions and thus produce a regular rational framing motion. 

\begin{theorem}
  \label{th:convex-hull}
  A rational spherical motion $\mathcal A(u,v)$ can be extended to a rational framing motion of a bounded regular rational curve if and only if the convex hull of the spherical trajectory of one point contains the origin.
\end{theorem}

\begin{proof}
Considering a suitable coordinate system we can assume that the spherical trajectory in the theorem statement is the trajectory of the point~$\qi$. Then the possibility to extend the spherical motion $\mathcal{A}(u,v)$ to the framing motion of a bounded rational curve is equivalent to the solvability of Problem~\ref{problem}, that is, the existence of a bounded positive rational function $L(u,v)$ such that, with $\mathbf{T}(u,v)$ defined as in \eqref{eq:T}, the integral \eqref{int} is rational.

  Now, if there is a rational framing motion of a bounded rational curve with spherical motion component given by $\mathcal{A}(u,v)$ then, as mentioned in Corollary~\ref{cor:zeroInt}, \eqref{zeroInt} holds true for a suitable bounded positive rational function $L(u,v)$. But this function can be interpreted as a weight function showing that the origin $\mathbf{0}$ lies indeed in the convex hull of $\mathbf{T}(u,v)$ -- the trajectory of the point~$\qi$.

  On the other hand, the convex hull property is equivalent to the existence of a positive continuous function $\lambda(u,v)\colon \mathbb S_1 \to \R$ so that
  \begin{equation}
    \int_{\mathbb S_1}\lambda(u,v)\mathbf T(u,v) \dif\mathbb S_1=\mathbf 0.
  \end{equation}
 We now choose arbitrarily a homogeneous quadratic polynomial $Q(u,v)=au^2+buv+cv^2$ that is irreducible over the real numbers. Hence, it has a pair of complex conjugate roots. Denote by $\mathcal C_0(\mathbb S_1)$ the algebra of all continuous real functions on the unit circle and, for fixed $\mathcal A(u,v)$ and $Q(u,v)$, consider its sub-algebra
  \begin{equation*}
    \mathcal R= \Bigl\{P(u,v) \frac{(\mathcal A(u,v) \mathcal A^*(u,v))^m (u^2+v^2)^n}{Q(u,v)^k} \mid P(u,v)\in\mathbb R[u,v];\ m,n,k \in{\mathbb N} \Bigr\}.
  \end{equation*}
  This algebra separates points and does not vanish at any point meaning that for any two points $(u_0,v_0)$, $(u_1,v_1)$ there exists $G(u,v)\in \mathcal R$ so that $G(u_0,v_0)\neq G(u_0,v_0)$ and for any point $(u_0,v_0)$ there exists $G(u,v)\in \mathcal R$ so that $G(u_0,v_0)\neq 0$. These are the two assumptions of the Stone-Weierstrass theorem \cite{deBranges1959} whence consequently $\mathcal R$ is dense in $\mathcal C_0(\mathbb S_1)$ in the topology of uniform convergence.

  Define the linear map $b\colon \mathcal C_0(\mathbb S_1)\to \mathbb R^3$
  \begin{equation}
    b(G) = \int_{\mathbb S_1}G(u,v)\mathbf T(u,v) \dif\mathbb S_1.
  \end{equation}
  By Lemma \ref{lemBanach} $(\ker b)\cap \mathcal R$ is dense in $(\ker b)$. Since $\lambda(u,v)\in (\ker b)$ it can be arbitrarily well approximated (in the sense of uniform convergence) by a function from $(\ker b)\cap \mathcal R$. Moreover $\lambda(u,v)>\epsilon>0$ and consequently there exists a function in $L(u,v)\in (\ker b)\cap \mathcal R$ which is positive.

  Now, let us apply Lemma~\ref{zeroRac} three times setting $F(u,v)$ to be successively each of the three components of
  \begin{equation}
    L(u,v)\mathbf T(u,v).
  \end{equation}
  Consequently $L(u,v)$, via the integral \eqref{int}, gives rise to a suitable parametric curve as trajectory of the origin for the extended framing motion.
\end{proof}

\begin{remark}
  \label{rem:PFD}
  The proof of Theorem~\ref{th:convex-hull} is not constructive but it provides an algebra of functions that is guaranteed to contain solutions. This algebra depends on the arbitrary choice of a quadratic irreducible polynomial $Q(u,v)$ and a suitably bounded regular PH curve $\mathbf{r}(u,v)$ will have an integer power of $Q(u,v)$ in its denominator. Call this an \emph{elementary rational PH curve.} Because the set of solution curves is closed under positive linear combinations, products of powers of arbitrary quadratic irreducible polynomials may arise in the denominator. Conversely, by the partial fraction decomposition of PH curves \cite{schroecker23}, any solution curve can be written as sum of elementary, but not necessarily regular, rational PH curves. We will return to this observation when computing bounded regular rational PH curves in Section~\ref{sec:construction}.
\end{remark}

\begin{figure}
  \centering
  \includegraphics{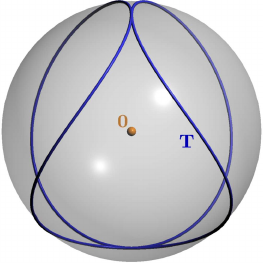}
  \caption{Spherical trajectory (tangent indicatrix) to Example~\ref{ex:T}}
  \label{fig:3}
\end{figure}

\begin{example}
  \label{ex:T}
  We consider the quaternion polynomial
  \begin{equation*}
    \mathcal{A}(u,v) = u^3 + (2\qj+\qk)u^2v - (1+2\qi)uv^2 - \qk v^3.
  \end{equation*}
  The spherical curve
  \begin{multline*}
    \mathbf{T}(u,v) =
    \frac{\mathcal{A}(t)\qi\Cj{\mathcal{A}}(t)}
    {\mathcal{A}(t)\Cj{\mathcal{A}}(t)} \\
    = \frac{(u^6-7u^4v^2+7u^2v^4-v^6)\qi+(2u^5v-12u^3v^3+2uv^5)\qj+(-4u^5v+4uv^5)\qk}{(u^2+v^2)^3}
  \end{multline*}
  together with the unit sphere and the origin is depicted in Figure~\ref{fig:3}. Its convex hull obviously contains the origin so that closed, regular, rational PH curves with this tangent indicatrix do exist. We will compute some of them in later examples in Section~\ref{sec:examples}.
\end{example}

\section{Computation of Regular Bounded Rational PH curves}
\label{sec:construction}

With the theoretical basis for closed PH curves established in Section \ref{sec:closed}, we now address computational aspects. Theorem~\ref{th:convex-hull} ensures the existence of bounded regular rational PH curves with a prescribed tangent indicatrix (containing the origin in its convex hull) and also for the existence of their Euler-Rodriguez framing motions. Because its proofs appeals to the Stone-Weierstrass theorem in a non-constructive way, it cannot be directly used for actual computations. In this section we complement Theorem~\ref{th:convex-hull} by a framework for the computation of bounded regular rational PH curves.

\subsection{Constructions over $\mathbb R$ and Residua Conditions}

In the theoretical context of Section~\ref{sec:closed}, there is no natural way to designate a special or exceptional point on a  bounded regular rational PH curve. In practice, however, one might want to construct these curves in a parametrization $\mathbf r(t)$ over the real line so that they only close in the limit at an exceptional point $\mathbf{r}(\infty) \coloneqq \lim_{t \to \pm \infty} \mathbf r(t)$. This point of view is not necessary but in line with basically all existing literature on PH curves. We will therefore transfer the results of the preceding section to this setting. The following corollary is the corresponding re-formulation of Theorem~\ref{th:all-CRRPH-curves} for an inhomogeneous parametrization.

\begin{theorem}
  \label{th:all-CRRPH-curves}
  All bounded regular rational PH curves can be obtained via the formula 
  \begin{equation}
    \label{int1}
    \mathbf r(t)=\int\ell(t)\frac{\mathcal A(t) \mathbf i \mathcal A^*(t)}{\mathcal A(t) \mathcal{A}^*(t)} \dif t
  \end{equation}
  where $\mathcal A(t) \in \H[t]$ is a quaternion polynomial and $\ell(t) \in \R(t)$ is a positive rational function of degree $-2$ such that the integral in \eqref{int1} is rational.
\end{theorem}

\begin{proof}
  In view of Theorem~\ref{th:all-CRRPH-curves} and its proof, all we need to show is that $\deg\ell(t) = -2$ is necessary and sufficient for boundedness. By virtue of \eqref{dS1} the integral \eqref{int1} transforms to
  \begin{equation*}
    \int_{t_0}^{t_1}
    \frac{2}{t^2+1}
    L(t)
    \mathbf T(t)
    \dif t.
  \end{equation*}
  We see that $\ell(t) = \frac{2}{t^2+1}L(t)$ is indeed of degree $-2$ if and only if $L(t)$ is of degree~$0$.
\end{proof}

For actual computations we prefer a slightly different formulation of \eqref{int1} where the polynomial $\mathcal{A}(t)\Cj{\mathcal{A}}(t)$ is absorbed into the denominator of $\ell(t)$:

\begin{corollary}
  \label{cor:all-CRRPH-curves}
  All bounded regular rational PH curves can be obtained via the formula 
  \begin{equation}
    \label{int2}
    \mathbf r(t)=\int\varkappa(t) \mathcal A(t) \mathbf i \Cj{\mathcal A}(t)
  \end{equation}
  where $\mathcal A(t) \in \H[t]$ is a quaternion polynomial and $\varkappa(t) \in \R(t)$ is a positive rational function of degree $-2\deg\mathcal{A}(t)-2$ such that the integral in \eqref{int2} is rational.
\end{corollary}

\begin{remark}
  It is possible to assume that the polynomial $\mathcal{A}$ in Equation~\ref{int2} is $\qi$-reduced. In fact, this should be assumed to not unnecessarily increase the degree of $\mathcal{A}(t)$ and hence also that of the numerator and denominator polynomials of the rational function~$\varkappa(t)$.
\end{remark}

Now, given $\mathcal{A}(t)$, the question remains how to select $\varkappa(t)$ in \eqref{int2} subject to two conditions:
\begin{enumerate}
\item The integral \eqref{int2} is rational and
\item $\varkappa(t)$ is positive.
\end{enumerate}
The first condition can be formulated in terms of linear constraints on the coefficients of the numerator of $\varkappa(t)$, the second condition has a formulation in terms of positive definiteness of a matrix. Both conditions can be combined easily to form a semidefinite optimization problem \cite[Section~B.3]{netzer23}.

Let us consider rationality of the integral \eqref{int2} first. We set
\begin{equation*}
  \varkappa(t) \coloneqq \frac{\mu(t)}{\alpha(t)}
\end{equation*}
with some prescribed denominator polynomial $\alpha(t) = Q_1(t)^{m_1} Q_2(t)^{m_2} \cdots Q_a(t)^{m_a}$ where $Q_i(t)$ is a real, irreducible, quadratic polynomial, positive on $\mathbb R$, and $m_i \in \mathbb{N}$ is the multiplicity of its conjugate complex zeros $z_i$, $\overline{z}_i$ (c.f. Remark~\ref{rem:PFD}). The coefficients $\mu_0$, $\mu_1$, \ldots, $\mu_m$ of the numerator polynomial $\mu(t) = \sum_{i=0}^m \mu_it^i$ are considered to be yet undetermined. Note that $\alpha(t)$ is of even degree and so is $\mu(t)$, by virtue of Corollary~\ref{cor:all-CRRPH-curves},
\[m=2(m_1+\cdots +m_a-\deg\mathcal{A}(t)-1).\]

In order to ensure rationality of \eqref{int1}, it is necessary and sufficient that the residua $\varkappa(t)\mathcal{A}(t)\qi\Cj{\mathcal{A}}(t)$ at the zeros of $z_1$, $\overline{z}_1$, $z_2$, $\overline{z}_2$, \ldots, $z_a$, $\overline{z}_a$ of $Q_1(t)$, $Q_2(t)$, \ldots, $Q_a(t)$ vanish:
\begin{equation}
  \label{eq:zero-residue}
  \residue_{z_i} \varkappa(t)\mathcal{A}(t)\qi\mathcal{A}^*(t) =
  \residue_{\overline{z}_i} \varkappa(t)\mathcal{A}(t)\qi\mathcal{A}^*(t) = 0,
  \quad
  i = 1,2,\ldots,a.
\end{equation}
This gives a system of linear equations for the unknown coefficients of $\mu(t)$. The equations have complex coefficients because of $z_i$, $\overline{z}_i \in \mathbb{C}$ but $\mu(t)$ is supposed to be real. In our implementation this is ensured by considering real and imaginary parts of the conditions in \eqref{eq:zero-residue}. Hence, \eqref{eq:zero-residue} gives a total of $2a$ linear equations over the real numbers.

It follows from Theorem~\ref{th:convex-hull} that with the $a$ pairs of conjugate complex zeros of $\alpha(t)$ fixed, the system of equations resulting from \eqref{eq:zero-residue} will admit non-trivial solutions for sufficiently high degree $m$ of $\mu(t)$ or, equivalently, for sufficiently high multiplicities $m_1$, $m_2$, \ldots, $m_a$ of the zeros of $\alpha(t)$. Moreover, it follows from the general theory of rational PH curves (cf. \cite{schroecker23b}) that we will find non-polynomial solutions only if one of the multiplicities $m_1$, $m_2$, \ldots, $m_a$ is sufficiently large (at least three in generic cases). For lower multiplicities the polynomial $\mu(t)$ will be divisible by factors of the denominator and the result will be a polynomial curve which indeed can not be closed. Sufficiently high multiplicities $m_1$, $m_2$, \ldots, $m_a$ are therefore necessary. 

Let us illustrate this this with an example.

\begin{figure}
  \centering
  \includegraphics{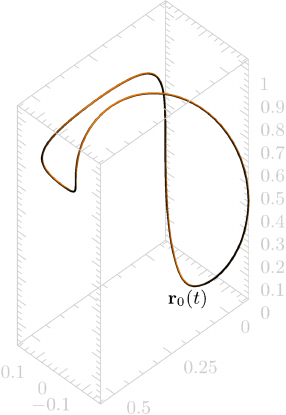}\quad\includegraphics{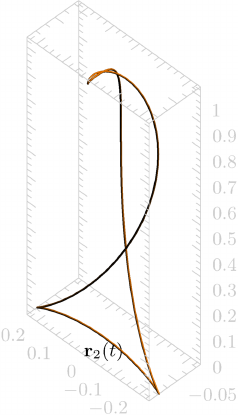}
  \caption{Closed, bounded rational PH curves $\mathbf{r}_0(t)$, $\mathbf{r}_2(t)$ of Example~\ref{ex:2}. Only $\mathbf{r}_0(t)$ is regular.}
  \label{fig:4}
\end{figure}

\begin{example}
  \label{ex:2}
  We consider the polynomial $\mathcal{A}(t) = t^3 + (2\qj + \qk)t^2 - (1 + 2\qi)t - \qk$. It is the de-homogenizaton of the polynomial of the same name of Example~\ref{ex:T}. Moreover, we write \[\varkappa(t) = \frac{\mu(t)}{\alpha(t)}\] where $\alpha = (4+t^2)^n$, $n \in \mathbb{N}$ and $\mu = \sum_{i=0}^m c_it^i$ with $m = 2(n-\deg\mathcal{A}(t)-1)$ and yet to be determined coefficients $c_0$, $c_1$, \ldots, $c_m$.

  In order for $m$ to be non-negative, we require at least $n \ge 4$. For $n = 4$ we have $m=0$ and a constant $\mu(t) = c_0$. The residua of $\varkappa(t) \mathcal{A}(t)\qi\Cj{\mathcal{A}}(t)$ at the zeros $z_1 = 2\ci$, $\overline{z}_1 = -2\ci$ of $\alpha(t)$ are complex scalar multiples of $c_0$. Thus, $\mu = 0$ and we only get trivial solutions.

  For $n = 5$, $m=2$ we have $\mu(t) = c_0 + c_1t+ c_2t^2$. The conditions for the vanishing of the residua are
  \[
      89c_0+7036c_2=0
      \quad\text{and}\quad
      c_1 = 0.
  \]
  This homogeneous system gives $\mu(t) = -89t^2 + 7036$ (or any multiple), a polynomial with two real roots. Hence the resulting closed rational PH curve is not regular but contains cusps. It is subsumed in our next try, $n = 6$ and $m=4$. The vanishing of the residua is equivalent to 
  \begin{equation*}
    11c_0 - 756 c_2 - 53264 c_4 = 0,\quad
    17 c_1 - 172 c_3 = 0,\quad
    41 c_1 + 436 c_3 = 0.
  \end{equation*}
  Note that these are also the vanishing conditions of the integral
  \eqref{intvanish} because we only have one pair of conjugate complex zeros in the denominator. We find
  \begin{equation}
    \label{linsol}
    c_1 = c_3 = 0
    \quad\text{and}\quad
    c_4 = \tfrac{11}{53264}c_0 - \tfrac{189}{13316}c_2.
  \end{equation}
  Indeed, for $c_0 = 1$, $c_2 = 0$, the integral \eqref{zeroInt} is rational:
  \begin{multline*}
    \mathbf{r}_0(t) = \frac{1}{798960(t^2+4)^5}
    \bigl(
    -165\qi t^9
    +(-165\qj+330\qk)t^8
    -1595\qi t^7\\
    +(-825\qj+2640\qk)t^6
    -168955\qi t^5
    +(-269675\qj+543090\qk)t^4\\
    +316255\qi t^3
    +(659090\qj+1086180\qk)t^2
    -199740\qi t
    +367480\qj+549360\qk
    \bigr).
  \end{multline*}
  For $c_0 = 0$, $c_2 = 1$, we obtain
  \begin{multline*}
    \mathbf{r}_2(t) = \frac{1}{798960(t^2+4)^5}
    \bigl(
      11340\qi t^9
      +(11340\qj-22680\qk)t^8
      -156700\qi t^7\\
      +(-342780\qj+617520\qk) t^6
      +256900\qi t^5
      +(230580\qj+2477640\qk)t^4\\
      -66580\qi t^3
      +(261420\qj+4555800\qk) t^2
      +209136\qj+3644640\qk
    \bigr).
  \end{multline*}
  The curves $\mathbf{r}_0(t)$ and $\mathbf{r}_2(t)$ are depicted in Figure~\ref{fig:4}. Apparently $\mathbf{r}_0(t)$ is regular while $\mathbf{r}_2(t)$ is not, an observation that can immediately be confirmed by looking at the respective numerator polynomials $\mu_0(t)$ and $\mu_2(t)$.
  
  The vector space of all solution curves is spanned by $\mathbf{r}_0(t)$, $\mathbf{r}_2(t)$, and the trivial polynomial solutions $\qi$, $\qj$, $\qk$. The latter only account for the integration constant and effect a translation of the solution curve. The set of all regular solutions is a convex cone in this vector space.
\end{example}

\subsection{Enforcing Speed Positivity}

Example~\ref{ex:2} confirms that it is really necessary to address the issue of curve regularity. In order to avoid singular points (cusps), we need to ensure that $\mu(t)$ is positive. A famous result of Hilbert \cite{hilbert88} guarantees that this is equivalent to $\mu(t)$ being the sum of squares of polynomials. Clearly, there exists a symmetric matrix $M$, not unique, of dimension $m/2 + 1$ where $m = \deg \mu(t)$ (note that $m$ is even) such that
\begin{equation}
  \label{eq:Gram-matrix}
  \mu(t) = (t^0, t^1, \ldots, t^{m/2}) \cdot M \cdot (t^0, t^1, \ldots, t^{m/2})^T.
\end{equation}
The polynomial $\mu(t)$ is the sum of squares of polynomials if and only if the matrix $M$ can be chosen positive definite \cite[Theorem~2.4]{choi94}. Thus, our problem can be re-formulated as follows: Determine $\mu(t)$ of degree $m$ as in
\eqref{eq:Gram-matrix} such that
\begin{itemize}
\item the coefficients of $\mu(t)$ or, equivalently, the coefficients of $M$
  fulfill a set of linear equations (coming from the zero residue constraints)
  and
\item $M$ is positive definite.
\end{itemize}
In this way, solutions to Problem~\ref{problem} are interpreted as points in the interior of the feasible region of a semidefinite program (SDP) \cite[Section~B.3]{netzer23} and we can rely on existing software for systematically computing them.

Denote by $x$ the vector of unknowns in the matrix $M$. SDP solves optimization
problems of type
\begin{equation*}
  \text{Minimize}\quad
  \langle c, x \rangle
  \quad
  \text{subject to}\quad
  M \succ 0
\end{equation*}
where $c$ is an arbitrary real vector of dimension equal to the dimension of $x$. In particular, the choice $c = 0$ yields a constant cost function so that the optimization problem is solved exactly by the interior feasible points. In principle, this is what we want but details need to be taken care of. We discuss this in the next section at hand of examples.

\section{Examples}
\label{sec:examples}

To begin with, we view an already discussed example and through the lens of our SDP formulation.

\begin{example}
  \label{ex:3}
  We continue Example~\ref{ex:2} and consider the case $n = 6$. Substituting \eqref{linsol} into $\mu(t)$, we
  may write $\mu(t) = (1,t,t^2) \cdot M \cdot (1,t,t^2)^T$ where $M = x_0M_0 +
  x_1M_1 + x_2M_2$ and
  \begin{equation*}
    M_0 =
    \begin{pmatrix}
      \frac{1512}{11} & 0 & 1 \\
      0 & 0 & 0 \\
      1 & 0 & 0
    \end{pmatrix},
    \quad
    M_1 =
    \begin{pmatrix}
      \frac{756}{11} & 0 & 0 \\
      0 & 1 & 0 \\
      0 & 0 & 0
    \end{pmatrix},
    \quad
    M_2 =
    \begin{pmatrix}
      \frac{53264}{11} & 0 & 0 \\
      0 & 0 & 0 \\
      0 & 0 & 1
    \end{pmatrix}.
  \end{equation*}
  The semidefinite program to be solved is
  \begin{equation}
    \label{eq:ex3}
    \text{Minimize}\quad
    0\quad
    \text{subject to}\quad
    x_0M_0 + x_1M_1 + x_2M_2 \succ 0.
  \end{equation}
  The point
  \begin{equation*}
    (x_0,x_1,x_2) = (-\tfrac{1}{100}, \tfrac{1}{50}, \tfrac{11}{53264})
  \end{equation*}
  leads to $\mu_0(t)$ and is feasible. Precisely the points of the shape
  \begin{equation*}
    (x_0,x_1,x_2) = (\chi, -2\chi, -\tfrac{189}{13316})
  \end{equation*}
  with $\chi \in \R$ arbitrary lead to $\mu_2(t)$. None of them is feasible. This confirms the visual appearance that the curve $\mathbf{r}_0(t)$ in Figure~\ref{fig:2} is free of cusps while $\mathbf{r}_2(t)$ is not. Numeric solution of \eqref{eq:ex3} using the Python library CVXPY \cite{diamond2016,agrawal2018} yields a solution curve that is visually indistinguishable from the regular solution of Figure~\ref{fig:2}.
\end{example}

\begin{remark}
  At this point, a few remarks on the numerics of the computation of bounded regular rational PH curves via semidefinite programming seem appropriate:
  \begin{itemize}
  \item Since SDP algorithms are inherently numerical, solutions polynomials $\mu(t)$ are only guaranteed to be ``almost'' positive. The effect of this numeric inaccuracy on the curves' visual appearance can be disastrous. Even if $\mu(t)$ is positive in a strict mathematical sense, the proximity to real zeros may result in curve regions of extremely high curvature. Ways to remedy this issue will be discussed later in Example~\ref{ex:3}.
  \item Since solutions are scale invariant, it is a good idea to fix one of the unknown coefficients of $\mu(t)$ prior to optimization. Nonetheless, the size of the resulting closed, regular, rational PH curves is hard to control. It will be necessary to scale them appropriately. The curves in all of our figures are scaled to approximately fit into a cube of edge length~$1$.
  \item Numeric solutions for $\mu(t)$ will result in almost vanishing residues of $\frac{\mu(t)}{\alpha(t)}\mathcal{A}(t)\qi\mathcal{A}^*(t)$. In other words, the integral curves are not rational. Pruning residues from a partial fraction decomposition of the integrand will resolve this issue but at the cost of not preserving the given tangent indicatrix. A better way is a mixed, symbolic-numeric approach where the coefficients of $\mu$ are substituted into a symbolic expression of the closed, regular, rational PH curve $\mathbf{r}(t)$ that already takes into account the zero residue conditions but not yet curve regularity.
  \end{itemize}
\end{remark}

It is not a priori clear that the choice $\alpha(t) = (4+t^2)^6$ in Examples~\ref{ex:2} and \ref{ex:3} will lead to solutions to Problem~\ref{problem}. However, the visual inspection of Example~\ref{ex:T} already confirmed that the origin is contained in the convex hull of the tangent indicatrix. Thus, Theorem~\ref{th:convex-hull} guarantees existence of solutions to $\alpha(t) = Q(t)^m$ where $Q(t) \in \R[t]$ is any quadratic irreducible polynomial and $m \in \mathbb N$ is \emph{sufficiently large.} We illustrate this with another example.

\begin{figure}
  \centering
  \includegraphics{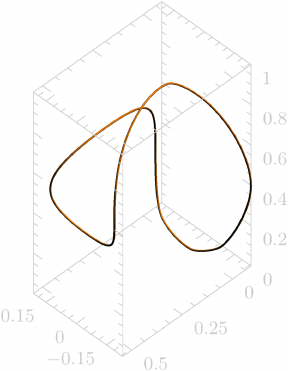}\hfill \includegraphics{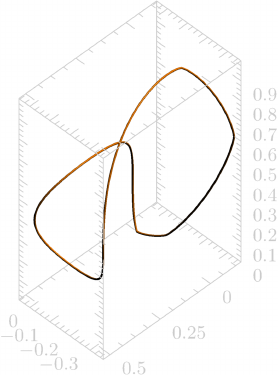}\hfill \includegraphics{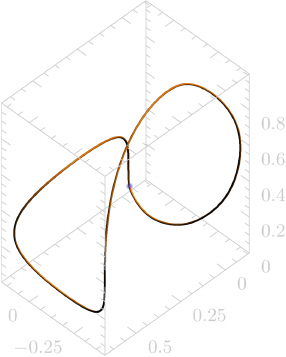}\newline
  \mbox{}\qquad\includegraphics{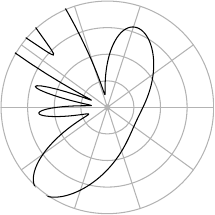}\hfill \includegraphics{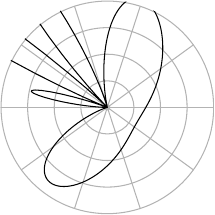}\hfill \includegraphics{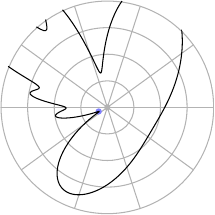}\qquad\mbox{}
  \caption{Solution curves to Example~\ref{ex:4}.}
  \label{fig:5}
\end{figure}

\begin{figure}
  \centering
  \includegraphics[width=6cm]{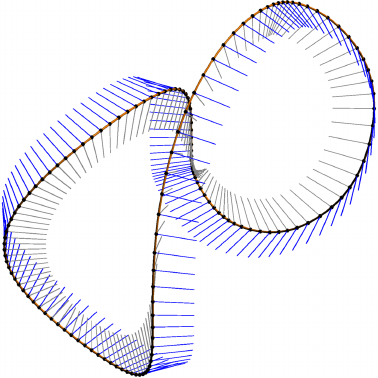}
  \caption{Closed regular rational framing motion.}
  \label{fig:6}
\end{figure}

\begin{example}
  \label{ex:4} We use the polynomial $\mathcal{A}(t)$ of Examples~\ref{ex:2} and \ref{ex:3} but select $\alpha(t) = (t^2-2t+5)^6$. A computation similar to that of Example~\ref{ex:2} produces a two parametric set of solution curves (up to translation) but all of them with cusps.
  
  Still, we know by Theorem~\ref{th:convex-hull} that a regular solution exists for $\alpha(t)$ equal to a sufficiently high power of $t^2-2t+2$. Indeed, for $\alpha(t) = (t^2-2t+2)^9$ semidefinite optimization suggests that regular solutions do exist. Two solution curves for different cost functions are displayed in Figure~\ref{fig:5}, left and center.
  
  The polar plot below each curve displays the length of the derivative vector as a function over the unit circle. It shows that the curve on the left is regular while this cannot be visually confirmed for the curve in the center. It is regular, at least numerically but this is not the point. Since it is too close to a singular curve, its visual appearance is unsatisfactory.
  
  To a certain extent, this is also true for the curve on the left. A natural way to improve visual appearance is to average several solution curves with positive weights. However, a study of polar plots suggest that all solution curves we compute are close to zero in the same regions over the unit circle so that no visual improvement is obtained. Increasing the degree of $\alpha(t)$, that is, setting $\alpha(t) = (t^2-2t+2)^{10}$ we obtain further degrees of freedom that provide the visually appealing curve in the right of Figure~\ref{fig:5}. Its polar graph indicates one potentially problematic point that is also indicated in the curve plot.

  Finally, Figure~\ref{fig:6} displays a closed regular rational framing motion to the curve of Figure~\ref{fig:5}, right. There are no singularities and the closed motion is infinitely differentiable.
\end{example}

Our use of SDP for computing solutions to Problem~\ref{problem} should not imply that this is the only possible approach. An alternative would be Cylindrical Algebraic Decomposition (CAD). This is a method to deal with sets described by polynomial equations and inequalities over the real numbers. It can be used to decide whether the set has interior points and also to determine interior sample points \cite{kauers11}. For the problem at hand, we need to augment the linear constraints on the coefficients of the matrix $M$ by strict inequalities that encode positive definiteness of $M$. One advantage of CAD is that the output will be exact if the input is symbolic, a situation often encountered in the context of PH curves. A disadvantage is that it is potentially very slow. In fact, we found that the general purpose implementation in the computer algebra system Mathematica is capable of producing only solutions for very simple examples.

The article \cite{ambrozie07} re-formulates the problem of finding a positive definite matrix whose entries are subject to linear constraints as a convex optimization problem. Existence of solutions is equivalent to existence of a minimizer of a certain convex function. Since analytic expressions for both the function and its gradient are available, the minimizer can be found numerically via gradient descent. This method computes solutions fast and efficiently but requires some care when evaluating the convex function or its gradient as both involve potentially rather large matrix exponentials. It could be used to find solutions to Problem~\ref{problem} but also to improve approximate solutions found by SDP. One disadvantage is that it only provides a single solution so that averaging of solutions is not possible.

\section{Conclusion}
\label{sec:conclusion}

We have presented a method to design space curves with closed framing motions which, in contrast to previous approaches \cite{FaroukiClosedC1,FaroukiClosedC2},  are infinitely differentiable. This necessitates that the origin of the moving frame moves along a bounded \emph{rational} PH curve. While rationality of this curve can easily be guaranteed by imposing zero-residue constraints, its regularity requires solution of a semidefinite program. Moreover, not every rational tangent indicatrix gives rise to regular solutions. A necessary and sufficient condition is that it contains the origin in its convex hull (Theorem~\ref{th:convex-hull}).

There are some open questions left for future research. First of all, there is high interest in low degree examples of rational PH curves giving rise to closed bounded rational framing motions. What is the minimal degree that can be achieved?

Moreover, given a finite collection $\mathbf{r}_1(t)$, $\mathbf{r}_2(t)$, \ldots, $\mathbf{r}_n(t)$, any linear combination $\mathbf{r}(t) = \sum_{i=1}^n \varrho_i\mathbf{r}_i(t)$ with positive coefficients $\varrho_1$, $\varrho_2$, \ldots, $\varrho_n$ is a regular solution again. It would be desirable to have tools that automatically produce positive linear combinations to optimize some curve fairness functional. Many typical fairness functionals are quadratic in $\varrho_1$, $\varrho_2$, \ldots, $\varrho_n$ \cite{nowacki96} and, in our context, allow to formulate the curve fairing problem as a quadratic program. Our hodograph based curve construction suggests to measure curve fairness by the integral of the squared length of the hodograph. Further investigation will be needed to determine its suitability in applied settings.

We have seen in Examples~\ref{ex:2} and \ref{ex:3} that the roots of the denominator polynomial $\alpha(t)$ can have significant influence on shape and degree of solution curves. It is unclear how to select $\alpha(t)$ in a way that ensures existence of ``good'' solutions of low degree. While we have discussed some aspects of solution quality, it is not even clear which solutions to consider ``good''. Low degree will certainly play a role as does curve fairness.

\bibliographystyle{elsarticle-num}

\end{document}